\documentclass[11pt]{amsart}

\newcommand{\N}{\mathbb{N}}
\newcommand{\Z}{\mathbb{Z}}
\newcommand{\R}{\mathbb{R}}

\newcommand{\T}{\mathbb{T}}
\def\di{\displaystyle}
\usepackage{amssymb}
\usepackage{latexsym}
\usepackage{xypic}
\usepackage{eufrak}
\usepackage{euscript}
\newtheorem{thm}{Theorem}
\newtheorem{defi}{Definition}

\def\cal{\mathcal}
\begin{document}
\baselineskip 6mm
\title{A $\lambda$-lemma for normally hyperbolic invariant manifolds}
\author{Jacky Cresson}
\address{Universit\'e de Franche-Comt\'e, \'Equipe de Math\'ematiques de Besan\c{c}on, CNRS-UMR 6623, 16 route de Gray, 25030
Besan\c{c}on cedex, France.} \email{cresson@math.univ-fcomte.fr}
\author{Stephen Wiggins}
\address{School of Mathematics, University of Bristol, University Walk, Bristol, BS8 1TW, United Kingdom.}
\email{S.Wiggins@bristol.ac.uk}

\maketitle

\begin{abstract}
Let $N$ be a smooth manifold and $f:N\rightarrow N$ be a $C^\ell$,
$\ell\geq 2$ diffeomorphism. Let $M$ be a normally hyperbolic
invariant manifold, not necessarily compact. We prove
an analogue of the $\lambda$-lemma in this case.
\end{abstract}

\tableofcontents

\section{Introduction}

In a recent paper, Richard Moeckel \cite{moe} developed a method
for proving the existence of drifting orbits on Cantor sets of
annuli. His result is related to the study of Arnold diffusion in
Hamiltonian systems \cite{ar}, and provides a way to overcome the
so called {\it gaps problems} for transition chains in Arnold's
original mechanism \cite{lo}. We refer to Lochak \cite{lo}
for a review of this problem.\\

The principal assumption of his work is the existence of a
symbolic dynamics for a compact normally hyperbolic invariant annulus.
His assumptions can be formulated as follows (\cite{moe},p.163):\\

Let $\Sigma$ denote the Cantor set of all bi-infinite sequences of
$0$'s and $1$'s, and $\sigma :\Sigma \rightarrow \Sigma$ be the
shift map defined for $s=(s_i )_{i\in \Z}$ as $\sigma (s)_i
=s_{i+1}$.

Let $N$ be a smooth manifold and $F:N\rightarrow N$ be a $C^\ell$
diffeomorphism. Let $M\subset N$ be a $C^\ell$ normally hyperbolic
invariant manifold for $F$ such that the stable and unstable
manifolds $W^{s,u} (M)$ intersect transversally in $N$. Then there
exists in a neighbourhood of $M$ an invariant hyperbolic set $\Lambda$ with the following properties:\\

i) $\Lambda \sim \Sigma \times M$.

ii) $\Lambda$ is a $C^{0,\ell}$ Cantor set of manifolds, i.e. each
leaf $M_s =\{ s\} \times M$, $s\in \Sigma$ is $C^\ell$ and depends
continuously on $s$ in the $C^\ell$ topology.

iii) some iterate $F^n \mid _{\Sigma}$ is conjugate to a skew product over the shift
$$\phi :\Sigma \times M \rightarrow \Sigma \times M,\ \phi (s,p)=(\sigma (s), g_s (p)),$$
where $g_s :M_s \rightarrow M_{\sigma (s)}$.\\

Moeckel \cite{moe} refers to previous work of Shilnikov \cite{sil}, Meyer and Sell \cite{ms} and Wiggins \cite{wi}. \\

Wiggins \cite{wi} proves an analogue of the Smale-Birkhoff theorem near a transversal homoclinic normally hyperbolic
invariant torus. However, this result can not be used to justify Moeckel's assumptions. Indeed, Wiggins's result is
based on:\\

$\bullet$ a particular normal form near the normally hyperbolic
invariant torus obtained by Shilnikov \cite{sil}.

$\bullet$ an annulus is not a compact boundaryless manifold, contrary to the torus.\\

Moreover, in most applications the compact annulus is obtained by
truncating a normally hyperbolic invariant cylinder, which is not compact.
However, non-compactness can be easily handled as it only provides
technical difficulties. This is not the case
when the normally hyperbolic invariant manifold has a boundary, which leads to technical as well as dynamical problems.\\

Other problems of importance deal with general compact
boundaryless invariant normally hyperbolic manifolds, as for
example normally hyperbolic invariant spheres. It has recently been shown
that normally hyperbolic invariant spheres are an important phase
space structure in Hamiltonian systems with three or more
degrees-of-freedom. Specific applications where they play a
central role are cosmology \cite{oast}, reaction dynamics
\cite{wwju}, \cite{wbw1}, and celestial mechanics \cite{wbw2}. \\

The proof of the Smale-Birkhoff theorem for normally hyperbolic
invariant tori by Wiggins \cite{wi} is based on a generalized
$\lambda$-lemma. This $\lambda$-lemma has been generalized by E.
Fontich and P. Martin \cite{fp} under more general assumptions
and $C^2$ regularity for the map.\\

There are three settings where a new type of $\lambda$-lemma would
be useful. What characterizes the difference in each case is the
geometrical structure of the normally hyperbolic invariant
manifold $M$.

\begin{itemize}

\item $M$ is non-compact, and can be characterized by a global
coordinate chart. This situation arises when we consider normally
hyperbolic invariant cylinders in Hamiltonian systems.
Non-compactness is dealt with by assuming uniform bounds on first
and second derivatives of certain functions (cf. \cite{GiLa}).

\item $M$ is compact with a boundary. This situation arises when
we truncate normally hyperbolic invariant cylinders to form
normally hyperbolic invariant annuli. The technical difficulty is
controlling the dynamics at the boundary.

\item $M$ is compact, but it cannot be described globally by a
single coordinate chart. This situation arises when we consider
normally hyperbolic invariant spheres.

\end{itemize}

In this paper, we prove a $\lambda$-lemma for normally hyperbolic
invariant manifolds, which are not necessarily compact. This
result allows us to prove a $\lambda$-lemma for a normally
hyperbolic annulus, i.e. for a normally hyperbolic compact
manifold with boundaries which is a subset of a non-compact
boundaryless normally hyperbolic manifold. We can also use the
same result to prove
a $\lambda$-lemma for compact invariant manifolds that cannot be described by a single coordinate chart.\\

The proof of the Smale-Birkhoff theorem as well as its applications for diffusion in Hamiltonian systems
will be studied in a forthcoming paper \cite{cw2}.

\section{A $\lambda$-lemma for normally hyperbolic invariant manifolds}

We first define the norms that we will use throughout this paper.
Essentially, we will only require two norms; one for vectors and
one for matrices. All our vectors can be viewed as elements of
$\R^n$ (for some appropriate $n$) and our matrices will consist of
real entries. As a vector norm we will use the sup norm on $\R^n$,
denoted by $\mid \cdot \mid$. Let ${\cal M}_{m\times n}$ denote
the set of $m\times n$ matrices over $\R$, $n\geq 1$, $m\geq 1$.
An element of ${\cal M}_{m\times n}$ has the form $A=(a_{i,j}
)_{i=1,\dots ,m,\ j=1,\dots ,n} \in {\cal M}_{m\times n}$. We
define  the norm of $A \in {\cal M}_{m\times n}$ by $\parallel
A\parallel =\di\sup_i \sum_j \mid a_{i,j} \mid$.

\subsection{Normally hyperbolic invariant manifolds}

Let $N$ be a $n$-dimensional smooth manifold, $n\geq 3$, and
$f:N\rightarrow N$ be a $C^\ell$ diffeomorphism, $r\geq 1$. Let
$M$ be a boundaryless $m$-dimensional submanifold (compact or non
compact) of class $C^\ell$ of $N$, $m< n$, invariant under $f$,
such that:\\

i) $M$ is normally hyperbolic, \\

ii) $M$ has a $m+n_s$-dimensional stable manifold $W^s (M)$ and a
$m+n_u$-dimensional unstable manifold $W^u (M)$, with
$m+n_s +n_u =n$.\\

Let $p\in N$, we denote by $Df_p$ the derivative of $f$ at $p$. Let $T_M N$ be
the tangent bundle of $N$ over $M$. As $M$ is normally hyperbolic, there exists
a $Df$-invariant splitting $T_M N=E^s \oplus E^u \oplus TM$ such that
$E^s \oplus TM$ is tangent to $W^s (M)$ at $M$ and $E^u \oplus TM$ is tangent to $W^u (M)$ at $M$. \\

\subsection{Normal form}
\label{sectionnormal}

We assume in the following that there exist a $C^\ell$ coordinate
systems $(s,u,x)\in \R^{n_s} \times \R^{n_u} \times M$ in a
neighbourhood $U$ of $M$ such that $f$ takes the form:
\begin{equation}
\label{normalform}
f(s,u,x)=(A_s (x)\, s , A_u (x)\, u ,g(x)) +r (s,u,x ) ,
\end{equation}
where $r$ is the remainder, $r =(r_s (s,u,x),r_u (s,u,x),r_x (s,u,x))$ with $A_s$, $A_u$, $g$ and $r$ satisfying\\

{\bf a) (invariance of $M$)} $r_s (0,0,x )=r_u (0,0,x)=r_x (0,0,x)=0$ for all $x\in M$.\\

As a consequence, the set $M$ is given in this coordinates system by
\begin{equation}
M=\{ (s,u,x)\in U \mid s=u=0 \} ,
\end{equation}
and $U$ can be chosen of the form
\begin{equation}
U=B_{\rho} \times M,
\end{equation}
with $\rho >0$ and $B_{\rho}$ is the open ball defined by $B_{\rho} =\{
(s,u)\in \R^{n_s} \times \R^{n_u} ;\
\mid (s,u)\mid < \rho \}$.\\

Let $\rho >0$, we denote by $B_{\rho}^s$ (resp. $B_{\rho}^u$) the open ball of
size $\rho$ in $\R^{n_s}$ (resp. $\R^{n_u}$) around
$0$.\\

As $M$ is normally hyperbolic, for $\rho >0$ sufficiently small,
the stable manifold theorem (see \cite{hps},\cite{wi2}) ensures
that the stable and unstable manifolds can be represented as
graphs, i.e. there exist two $C^\ell$ functions $G^s (s,x)$ and
$G^u (u,x)$ such that
\begin{equation}
\left .
\begin{array}{lll}
W^s (M)\cap U & = & \{ (s,x) \in B_{\rho}^s \times M \mid u=G^s (s,x) \} , \\
W^u (M)\cap U & = & \{ (u,x) \in B_{\rho}^u \times M \mid s=G^u (u,x) \} ,
\end{array}
\right .
\end{equation}
with
\begin{equation}
G^s (0,x)=0,\ \partial_s G^s (0,x)=0,\ \partial_x G^s (0,x)=0 ,
\end{equation}
and
\begin{equation}
G^u (0,x)=0,\ \partial_u G^u (0,x)=0,\ \partial_x G^u (0,x)=0 ,
\end{equation}
which reflect the tangency of $W^s (M)$ and $W^u (M)$ to $E^s$ and $E^u$ over $M$ respectively.\\

Using these functions, we can find a coordinate system for which the stable and unstable manifolds are ``straightened'', i.e.\\

{\bf b) (straightening of the stable manifold)} $r_u (s,0,x )=0$
for all $(s,0,x)\in U$,

{\bf c) (straightening of the unstable manifold)} $r_s (0,u,x)=0$ for all $(0,u,x)\in U$.\\

As a consequence, the stable and unstable manifolds of $M$ are
given by
\begin{equation}
\left .
\begin{array}{lll}
W^s (M) & = & \{ (s,u,x)\in U \mid u=0 \} ,\\
W^u (M) & = & \{ (s,u,x)\in U \mid s=0 \} .
\end{array}
\right .
\end{equation}
Indeed, following the classical work of Palis-deMelo \cite{pd}, the change of
variables
\begin{equation}
\Phi : \left .
\begin{array}{lll}
U & \longrightarrow & U \\
(s,u,x) & \longmapsto & (s-G^u (u,x) ,u-G^s (s,x) ,x ) ,
\end{array}
\right .
\end{equation}
realizes the straightening:\\

For all $P\in W^s (M)$, we denote by $\Phi(P)=(s',u',x')$. Then, $P\in W^s (M)$
if and only if $\Phi (P)=(s',0,x')$ and $P\in W^u (M)$ if and only if $\Phi (P)=(0,u',x')$.\\

{\bf d) (conjugacy on the stable and unstable manifold)} We assume
that $r_x (0,u,x)=0$ and $r_x (s,0,x)=0$ for all $s\in
B^s_{\rho}$,
$u\in B^u_{\rho}$ and $x\in M$.\\

In many examples of importance this condition is satisfied. It tells us that
the dynamics on the stable and unstable manifolds in the invariant manifold
direction is given by the dynamics on $M$. We refer to Graff \cite{gr} for such
an
example of rigidity in an analytic context.\\

{\bf e) (hyperbolicity)} $\parallel A_s(x)\parallel \leq \lambda <
1$, $\parallel A_u(x)^{-1}
\parallel \leq \lambda < 1$. \\

These results lead us to introduce the following definition of a normal form
for diffeomorphisms near a normally hyperbolic invariant manifold:

\begin{defi}[Normal form]
Let $N$ be a smooth manifold and $f$ a $C^\ell$ diffeomorphism of
$N$, $\ell\geq 2$. Let $M$ be a compact normally hyperbolic
invariant manifold of $f$. The diffeomorphism is said to be in
normal form if there exist a neighbourhood $U$ of $M$ and a
$C^\ell$ coordinate system on $U$ such that $f$ takes the form
(\ref{normalform}) and satisfies conditions a)-e).
\end{defi}

Standard results on normal form theory can be used to prove in some case that we have a diffeomorphism in normal form. We
refer to (\cite{bk},p.332) for a general normal form theorem. In particular, we derive such a normal form in a Hamiltonian
setting near a normally hyperbolic cylinder.\\

Moreover, general normal form results for normally hyperbolic manifolds already imply that our assumptions are general, at
least if we restrict the regularity assumption on the coordinates system to $C^1$. Indeed, we have the
following result due to M. Gidea and R. De Llave \cite{GiLa}:

\begin{thm}
Let $N$ be a smooth manifold and $f$ a $C^\ell$ diffeomorphism of
$N$, $\ell\geq 2$. Let $M$ be a normally hyperbolic invariant
manifold of $f$ (compact or non compact). There exists a
neighbourhood $U$ of $M$ and a $C^1$ coordinate system on $U$ such
that $f$ is in normal form.
\end{thm}

We refer to (\cite{GiLa},$\S$.5.1) for a proof.\\

Of course, such a result is not sufficient for our purposes as we
need some control on the second order derivatives of $g$ and $r$.
However, it proves that our assumptions are general. In the same
paper, M. Gidea and R. De Llave \cite{GiLa} proves that we can
take
$r_x =0$, i.e. that we have a decoupling between the center dynamics and the hyperbolic dynamics.\\

\subsection{The $\lambda$-lemma}

We have the following generalization of the toral $\lambda$-lemma of S. Wiggins \cite{wi}:

\begin{thm}[$\lambda$-lemma]
\label{lamlem} Let $N$ be a smooth manifold and $M$ be a $C^\ell$
submanifold of $N$, normally hyperbolic, invariant under a
$C^\ell$ diffeomorphism $f$, $\ell\geq 2$, in normal
form in a given neighbourhood $U$ of $M$ and such that\\

i) There exists $C>0$ such that
\begin{equation}
\label{secondorder}
\sup \left \{ \parallel \partial^2_{\sigma
,\sigma'} r_i (z) \parallel,\ z\in U , \sigma \in \{ s,u,x\},\,
\sigma' \in \{ u,x\}, \ i \in \{ s,x\} \right \} \leq C .
\end{equation}

ii) There exists $\tilde{C}>0$ such that
\begin{equation}
\sup \left \{ \parallel \partial^2_{\sigma ,x} g (z) \parallel ,\
z\in U ,\ \sigma \in \{ s,u,x\} \right \} \leq \tilde{C}.
\end{equation}

iii) There exists $D>0$ such that for all $x\in M$, $\parallel \partial_x A_s (x)\parallel \leq D$.\\

\noindent
Let $\Delta$ be an $m+n_u$ dimensional manifold
intesecting $W^s (M)$ transversally and let $\Delta_k =f^k \left(
\Delta \right) \cap U$ be the connected component of $f^k \left(
\Delta \right) \cap U$ intersecting $W^s \left( M \right)$. Then
for $\epsilon
>0$, there exists a positive integer $K$ such that for $k \geq K$, $\Delta_k$ is $C^1$ $\epsilon$-close to $W^u (M)$.
\end{thm}

The proof follows essentially the same line as in (\cite{wi},p.324-329) and is given in section \ref{prooflambda}.

\section{Hamiltonian systems and normally hyperbolic invariant cylinders and annuli}

Normally hyperbolic invariant annuli or cylinders are the basic
pieces of all geometric mechanisms for diffusion in Hamiltonian
systems. This may seem to be an unusual statement in light of the
fact that the classical ``transition chain'' is a series of
heteroclinic connections of stable and unstable manifolds of
``nearby'' lower dimensional tori. However, these lower
dimensional tori are contained in normally hyperbolic invariant
annuli and cylinders which have their own stable and unstable
manifolds (which, in turn, contain the stable and unstable
manifolds of the lower dimensional tori used to construct Arnold's
transition chains). The importance of normally hyperbolic
invariant annuli or cylinders can be clearly seen in the papers of
Z. Xia \cite{xia}, R. Moeckel \cite{moe2}, and A. Delshams, R. De
Llave and T. Seara
\cite{dls} where normally hyperbolic annuli are a fundamental tool.\\

In this section, we prove a $\lambda$-lemma for Hamiltonian
systems possessing, in a fixed energy manifold, a normally
hyperbolic manifold of the form $\T \times I$, where $I$ is a
given compact interval, which belong to a non-compact boundaryless
invariant manifold $\T \times \R$. We will see that the fact that
a $\lambda$-lemma can be proven in this case is related to the
fact that the boundaries are partially hyperbolic invariant tori
for which an analogue of the $\lambda$-lemma for which a different
type of $\lambda$-lemma has already been proven. Moreover, we will
also construct a class of three degree of freedom Hamiltonian
systems which satisfy our assumptions.

\subsection{Main result}
\label{mainham}

All of our results will be stated in terms of discrete time
systems, or maps. However, many of the applications we have in
mind will be for continuous time Hamiltonian systems. Our results
will apply in this setting by considering an appropriate
Poincar\'e map for the continuous time system. It is important to
keep this reduction from continuous time Hamiltonian system to
discrete time Poincar\'e map firmly in mind from the point of view
of considering the dimensions of the relevant invariant manifolds
and transversal intersection in the two systems. Finally, we note
that even though the specific applications we consider here are
for Hamiltonian systems, as is true of most hyperbolic phenomena,
a Hamiltonian structure is not generally required for their
validity. We now describe the setting and our hypotheses for the
applications of interest.\\

Let $H$ be a $C^r$ Hamiltonian  $H(I,\theta )$, $(I,\theta)\in
\R^3 \times \T^3$, $r\geq 3$. The Hamiltonian defines a $C^{r-1}$
Hamiltonian vector field for which we make the following
assumptions: \\

i) the Hamiltonian vector field possesses an invariant normally
hyperbolic manifold $\bar{\Lambda}$ which is diffeomorphic to
$\T^2 \times \R^2$, with $5$-dimensional stable and unstable
manifold $W^s (\bar{\Lambda} )$ and $W^u
(\bar{\Lambda} )$ respectively (note: these invariant manifolds are not isoenergetic).\\

ii) There exists a Poincar\'e section in a tubular neighbourhood
of $\bar{\Lambda}$ and a $C^2$ coordinate system of the form $(v,w
,s,u) \in S^1 \times \R \times \R \times \R$ such that the
Poincar\'e map is in normal form near $\bar{\Lambda}$, i.e. it
takes the form

\begin{equation}
\label{normalform2}
f(s,u,v,w)=(A_s (z)\, s , A_u (z)\, u ,g(v) ,w ) +r (s,u,v,w ) ,
\end{equation}

\noindent
 where $r$ is the remainder, $r =(r_s (s,u,v,w),r_u
(s,u,v,w),r_x (s,u,v,w))$ with $A_s$, $A_u$, $g$ and $r$
satisfying assumptions a)-e) of section \ref{sectionnormal}. We
denote the intersection of $\bar{\Lambda}$, $W^s (\bar{\Lambda})$,
and $W^u(\bar{\Lambda})$ with the Poincar\'e section by $\Lambda$,
$W^s (\Lambda)$, and $W^u(\Lambda)$, respectively.  The
$4$-dimensional Poincar\'e section is chosen such that the
intersection of these manifolds with the Poincar\'e section is
isoenergetic, and $\Lambda$ is $2$-dimensional, $W^s (\Lambda)$ is
$3$-dimensional, and $W^u(\Lambda)$ is  $3$-dimensional.

iii) There exist two circles $C_0$ and $C_1$ belonging to
$\Lambda$ and invariant under $f$, possessing $2$ dimensional
stable and unstable manifolds.

iv) We assume that a $\lambda$-lemma is valid for the partially hyperbolic invariant circles $C_0$ and $C_1$.\\

Assumptions iv) can be made more precise as there already exists
many versions of the $\lambda$-lemma for partially
hyperbolic tori. We refer in particular to \cite{cr1},\cite{cr3},\cite{ma} and \cite{fp}, which is the most general.\\

We denote by $A$ the invariant normally hyperbolic annulus whose
boundaries are $C_0$ and $C_1$.  $W^s (A)\subset W^s (\Lambda )$
is the stable of $A$ and  $W^u (A)\subset W^u (\Lambda )$ is the
unstable manifold of $A$ respectively. We define the boundary of
$W^s (A)$ to be $W^s(C_0)$ and $W^s (C_1)$ and the boundary of
$W^u (A)$ to be $W^ u (C_0)$ and $W^u (C_1)$. The important point
here is that even though $A$ has a boundary, it is still invariant
with respect to both directions of time. This is because its
boundary is an invariant manifold. Similarly, $W^s (A)$ and $W^u
(A)$ are also invariant manifolds (cf. with inflowing and
outflowing invariant manifolds with boundary described in
\cite{wi2}). \\

Our main result is the following:

\begin{thm}[$\lambda$-lemma for normally hyperbolic annuli]
Let $H$ denote the Hamiltonian for a three degree of freedom
Hamiltonian system satisfying assumptions i)-iv). Let $\Delta$ be
a $3$ dimensional manifold intersecting $W^s (\Lambda )$
transversally. We assume that there exists a subset
$\tilde{\Delta}\subset \Delta$, such that $\tilde{\Delta}$
intersects $W^s (A)$ transversally, and such that the boundaries
$\partial \tilde{\Delta}_0$ and $\partial \tilde{\Delta}_1$
intersect transversally the stable manifolds of $C_0$ and $C_1$
respectively. Then, for all $\epsilon >0$, there exists a positive
integer $K$ such that for all $k\geq K$, $f^k (\Delta )$ is $C^1$
$\epsilon$-close to $W^u (A)$. \label{lamlem3dof}
\end{thm}

The proof follows from our previous Theorem \ref{lamlem} for the
noncompact case and the $\lambda$-lemma for partially hyperbolic
tori to control the boundaries of $\tilde{\Delta}$. In considering
the proof of Theorem \ref{lamlem} one sees that the difficulty
arising for invariant manifolds is that iterates of points may
leave the manifold by crossing the boundary. This is dealt with
here by choosing the boundary of the manifold to also be (lower
dimensional) invariant manifold(s).

\begin{proof} Since the boundary of $W^s (A)$ is invariant, we know that
under iteration by $f$ $\tilde{\Lambda} \cap W^s
(A)$ is always contained in $W^s(A)$. Hence we can use the
 $\lambda$-lemma proven in Theorem \ref{lamlem} for
non-compact normally hyperbolic invariant manifolds to conclude
that for all $\epsilon >0$, there exists a positive integer $K$
such that for all $k\geq K$,   $f^k (\tilde{\Delta} )$ is  $C^1$
$\epsilon$-close to $W^u (A)$.

We also need to show that the boundaries of $\tilde{\Delta}$
correctly accumulate on the boundaries of $W^u (A)$ which are
formed by $W^u (C_0 )$ and $W^u (C_1 )$. This follows by applying
the $\lambda$-lemma for partially hyperbolic tori of Fontich and
Martin \cite{fp} to $C_0$ and $C_1$. It follows from this lemma
that for all $\epsilon >0$, there exists a positive integer $K'$
such that for all $k\geq K'$, $f^k (\partial \tilde{\Delta}_i )$
is  $C^1$ $\epsilon$-close to $W^u (C_i )$ for $i=0,1$. This
concludes the proof.
\end{proof}

We remark that if one is in a neighborhood of $\Lambda$ where the
manifolds are ``straightened'' it is likely that one can construct
$\tilde{\Delta}$ as a graph over $W^ u(A)$ such that it intersects
$W^s (A)$ transversally, and such that the boundaries $\partial
\tilde{\Delta}_0$ and $\partial \tilde{\Delta}_1$ intersect
transversally the stable manifolds of $C_0$ and $C_1$,
respectively, \\

In the following section we prove that our assumptions are satisfied in a large class of near integrable Hamiltonian
systems.

\subsection{An example}

In order to construct a normally hyperbolic cylinder for a near
integrable Hamiltonian system, we follow the same geometrical
set-up as in the seminal paper of Arnold \cite{ar}. We emphasize
that the results described above are not perturbative in nature,
but we have in mind near-integrable systems as a setting where our
results provide a key ingredient for proving the existence of
symbolic dynamics an instability mechanism.

\subsubsection{Variation around the Arnold example}

We consider the near integrable Hamiltonian system

\begin{equation}
H_{\epsilon ,\mu} (p,I,J,q,\theta ,\phi)={1\over 2} p^2 +  {1\over
2} I^2 +J+\epsilon (\cos q -1) +\epsilon f(\theta ,\phi ) +\mu
(\sin q)^{\alpha (\nu ,\sigma )} g(\theta ,\phi ), \label{modham}
\end{equation}

\noindent
 where as usual, $(I,J,\theta ,\phi )\in \R \times \R
\times \T \times \T$ are action-angle variables, $(p,q)\in \R
\times \T$, $0<\epsilon <<1$ is a small parameter and $\mu$ is
such that $0<\mu <<\epsilon$, and $f$ and $g$ are two given smooth
functions, and $\nu \in \N$ is a parameter controlling the order
of contact between $H_{\epsilon ,\mu}$ and $H_{\epsilon ,0}$ via
the function $\alpha (\nu ,\sigma )= 2\left [ \di {\log \nu \over
4\sigma } +1 \right ]$ introduced in (\cite{lm}, equation (2.5)),
$\sigma >0$ and $\nu \geq \nu_{\sigma}$ where $\nu_{\sigma}$ is
the smaller positive integer such that $\alpha (\nu )=2$ and
$\alpha (\nu ,\sigma )\geq 2$ for $\nu \geq \nu_{\sigma}$.\\

For $\epsilon=\mu=0$ the system is completely integrable and the set
\begin{equation}
\bar{\Lambda }=\{ (p,I,J,q,\theta ,\phi)\in \R \times \R \times \R \times \T\times \T \times \T \mid p=q=0 \} ,
\end{equation}
is invariant under the flow and normally hyperbolic. The dynamics on $\bar{\Lambda}$ is given by the completely integrable
Hamiltonian system
\begin{equation}
H(I,\theta )= {1\over 2} I^2 +J.
\end{equation}
As a consequence, the set $\bar{\Lambda}$ is foliated by invariant $2$-tori.\\

For $\mu=0$ and $\epsilon\not= 0$, the set $\bar{\Lambda}$
persists, but the dynamics on $\bar{\Lambda}$ is no longer
integrable. In particular, it is not foliated by invariant two
tori. However, the KAM theorem applies and we have a Cantor set of
invariant two tori whose measure tends to the full measure when
$\epsilon$ goes to zero. As a consequence, we are in a situation
where the ``large gaps'' problem arises (\cite{lo} ), contrary to
the well known example of Arnold \cite{ar} where all the foliation
by
invariant two tori is preserved under the perturbation.\\

Let $h$ be a given real number. We denote by ${\cal H}_{\epsilon
,0} =H_{\epsilon ,0}^{-1} (h)$ the energy manifold. There exists a
global cross-section to the flow denoted by $S$ and defined by

\begin{equation}
S=\{ (p,I,J,q,\theta ,\phi)\in \R \times \R \times \R \times \T\times \T \times \T \mid \phi =0 \} .
\end{equation}

\noindent We denote by $\Lambda$ the intersection of
$\bar{\Lambda}$ with $S \cap {\cal H}_{\epsilon ,0}$. We can find
a symplectic analytic coordinate system on $S\cap {\cal
H}_{\epsilon ,0}$, denoted by $(x,y,s,u)\in \T \times \R  \times
\R \times \R$ such that the set $\Lambda$ is defined by

\begin{equation}
{\Lambda} =\{ (x,y,s,u)\in \T \times \R \times \R \times \R ;\
s=u=0\, \} .
\end{equation}

\noindent
Geometrically, $\Lambda$ is a cylinder. This cylinder is a normally hyperbolic boundaryless manifold, but not compact.\\

We can introduce the compact counterpart, which is a normally hyperbolic annulus, but now with boundaries. Let $T_0$ and
$T_1$ be two invariant $2$-dimensional partially hyperbolic tori belonging to $\bar{\Lambda}$. We denote by $C_0$ and $C_1$
the intersection of $T_0$ and $T_1$ with $S\cap {\cal H}_{\epsilon ,0}$. The invariant circles $C_0$ and
$C_1$ are defined as
\begin{equation}
C_i =\{ (x,y,s,u) \in\T\times \R\times \R\times \R ;\ y=y_i ,\ s=u=0 \, \},
\end{equation}
for $y_i \in \R$ well chosen, $i=0,1$. We assume in the following that $y_0 < y_1$.

The compact counterpart of $\Lambda$ is then defined as
\begin{equation}
A =\{ (x,y,s,u)\in \T \times \R \times \R\times \R;\ y_0 \leq y\leq y_1  ,\ s=u=0\, \} .
\end{equation}

Let $P$ be the Poincar\'e first return map associated to $S\cap {\cal H}_{\epsilon ,0}$. The dynamics on $A$ is given by an
$\epsilon$ perturbation of an analytic twist map, i.e. that $P\mid_A$ is defined by
\begin{equation}
P\mid_A (x,y)=(x+\omega (y) ,y ) +\epsilon r(x,y),
\end{equation}
where $\omega' (y) \not= 0$ and $r(x,y_i )=0$ for $i=0,1$.\\

When $\epsilon\not=0$ and $\mu\not=0$, then the set $\Lambda$
persists since the perturbation vanishes on $\Lambda$.
However, the stable and unstable manifold of $\Lambda$ intersect transversally for a well chosen perturbation $g$.\\

We then obtain an example of a three degrees of freedom Hamiltonian system satisfying the geometrical assumptions i), iii)
and iv) of section \ref{mainham}. We prove in the following section that the analytic assumption ii) is satisfied.

\subsubsection{Poincar\'e section and normal form}

The main problem is to prove that the diffeomorphism of the
cross-section to the flow defined in a neighbourhood of $\Lambda$
is in normal form with respect to definition \ref{normalform}.\\

The basic theorem which we use to obtain a smooth normal form near
a normally hyperbolic manifold is  a generalized version of the
{\it Sternberg linearization theorem}. In the compact case, this
result has been obtain by Bronstein and Kopanskii (\cite{bk},
theorem 2.3,p.334). The non compact case has been proven by P.
Lochak and J-P. Marco \cite{lm}. This theorem, which can be stated
for flows or maps, ensures that we can obtain a conjugacy as
smooth as we want between $\phi_{\epsilon ,\mu}$ and
$\phi_{\epsilon ,0}$ by choosing the two flows with a sufficiently
high order of contact (see \cite{bk},p.334). A key remark in our
case is that the flows $\phi_{\epsilon, 0}$ and $\phi_{\epsilon
,\mu}$ generated by $H_{\epsilon ,0}$ and $H_{\epsilon ,\mu }$,
respectively, have  contact of order $\alpha (\nu ,\sigma )$ on
$\bar{\Lambda}$. Moreover, we can obtain an arbitrary order of
contact between the two flows by choosing the parameter $\nu$
sufficiently large. As a consequence, we will always be able to
realize the assumptions of the Sternberg linearization
theorem for normally hyperbolic manifolds and as a consequence, to obtain a normal form as smooth as we desire.\\

Before stating the normal form theorem, which is only a minor
modification of the result of Lochak-Marco (see
\cite{lm}, Theorem D), we introduce some notation:\\

As a general notation, for any manifold $M$, we denote by
$M_{\rho}$ a tubular neighbourhood of $M$ of radius $\rho$. We let
$f_{\nu}$ and $f_*$ denote the Poincar\'e maps defined in a
neighborhood of $\Lambda$ obtained from the flows generated by the
Hamiltonians $H_{\epsilon ,\mu}$ and $H_{\epsilon
,0}$, respectively.\\

Using the Sternberg linearization theorem for normally hyperbolic manifolds proved in \cite{lm}, we obtain the following
theorem:

\begin{thm}[Normal form]
\label{normalform3}
For $\nu_0$ large enought, there exist $\rho, \rho'$ with $0<\rho' <\rho$ such that for all $\nu\geq \nu_0$,
$f_{\nu} (\Lambda_{2\rho'} ) \subset
\Lambda_{\rho}$ and there exists a $C^k$ diffeomorphism $\phi_{\nu}$ ($k\geq 1$) satisfying $\Lambda_{\rho'} \subset
\phi_{\nu} (\Lambda_{2\rho'} )\subset \Lambda_{\rho}$ and $\phi_{\nu} \circ f_* =f_{\nu} \circ \phi_{\nu} $ on $\Lambda_{\rho'}$.
Moreover, there exists a constant $a$ ($0< a <1$)such that $\parallel \phi_{\nu}^{\pm 1} -Id \parallel_{C^k }
\leq a^{\alpha (\nu)}$, where $\parallel \cdot \parallel_{C^k}$ denotes the $C^k$ norm on $\Lambda_{\rho'}$.
\end{thm}

This theorem is a direct application of a result of Lochak-Marco
\cite{lm}, and we refer to their paper for more details (in
particular Theorem D).\\

Theorem \ref{normalform3} implies that the $\lambda$-lemma  proven in Theorem \ref{lamlem3dof} applies to the three degree-of-freedom Hamiltonian
systems defined by (\ref{modham}). \\

In this case one can, moreover, prove an analogue of the
Smale-Birkhoff theorem using the fact that symbolic dynamics is
stable under small $C^1$ perturbations and taking for the
Poincar\'e map near the normally hyperbolic invariant annulus the
linear mapping. A complete study of this problem will be done in
\cite{cw2} in the context of homoclinic normally hyperbolic
invariant manifolds.

\section{Proof of the $\lambda$-lemma}
\label{prooflambda}

\subsection{Preliminaries}

In this section we develop the set-up for the proof of the
$\lambda$-lemma. First we discuss some useful consequences of the normal form assumptions for the diffeomorphism $f$.\\

$\bullet$ We denote by $Dr (p)$ the differential of $r$ at point
$p\in U$. By invariance of $M$ (assumption a) in section
\ref{sectionnormal}) we have $Dr(p)=0$ for all $p\in M$. Let
$0<k<1$, since $r$ is a $C^1$ function, we have for $U$
sufficiently small and for all $p\in U$,

\begin{equation}
\parallel Dr (p) \parallel \leq k .
\label{k_est}
\end{equation}

\noindent
 This implies in particular that for all $p\in U$, the
partial derivatives $\partial_i r_j (p)$, $i,j\in \{ s,u,x\}$
satisfy $\parallel \partial_i r_j (p) \parallel \leq k$.\\

We take $U$ sufficiently small in order to have the following inequalities satisfied for $k$:\\

\begin{eqnarray}&&  0<\lambda +k <1, \label{lk1} \\
&& \lambda^{-1} -k >1. \label{lk2}
\end{eqnarray} \\

$\bullet$  The straightening conditions b) and c) imply that:
\begin{equation}
\label{implibc} \left .
\begin{array}{l}
\forall p\in W^s (M)\cap U,\ \ \ \partial_s r_u (p)=\partial_x r_u (p)=0 ,\\
\forall p\in W^u (M)\cap U,\ \ \ \partial_u r_s (p)=\partial_x r_s (p)=0.
\end{array}
\right .
\end{equation}

$\bullet$ The conjugacy assumption d) implies that:
\begin{equation}
\label{implid}
\left .
\begin{array}{l}
\forall p\in W^s (M)\cap U,\ \ \ \partial_s r_x (p)=\partial_x r_x (p)=0 ,\\
\forall p\in W^u (M)\cap U,\ \ \ \partial_u r_x (p)=\partial_x r_x (p)=0.
\end{array}
\right .
\end{equation}

\subsection{Notation}

Let $v_0$ be a unit vector in the tangent bundle to $\Delta_n$. We
denote the components of $v_0$ as $v_0 =(v_0^s ,v_0^u ,v_0^x )$
where $v_0^s \in E^s$, $v_0^u \in E^u$ and $v_0^x \in TM$. We
denote by $v_n$ the iterate of $v_0$ under $Df^n$, i.e. $v_n =Df^n
(v_0 )$ and $v_n =(v_n^s ,v_n^u ,v_n^x )$. Let $p=(s,u,x)$, we
denote by
$p_n =f^n (p)$ the $n$-th iterate of $p$ and $p_n =(s_n ,u_n ,x_n )$ its components. We identify $p_0$ and $p$.\\

Let $p\in U$, then using (\ref{normalform}) we have

\begin{equation}
Df_p = \left (
\begin{array}{lll}
A_s (x)+\partial_s r_s (p)& \partial_u r_s (p)& \partial_x r_s (p)+\partial_x A_s (x)\, s \\
\partial_s r_u & A_u (x ) +\partial_u r_u (p)& \partial_x r_u (p)+\partial_x A_u (x)\, u \\
\partial_s r_x (p)& \partial_u r_x (p)& \partial_x g (x)+\partial_x r_x (p)
\end{array}
\right ). \label{jacob1}
\end{equation}

We are going to prove that for all $v_0 \in T\Delta_n$, and for
$n$ sufficiently large, we have

\begin{equation}
\label{straight}
\lim_{n\rightarrow +\infty} \sup \left \{ I_n^x,
I_n^s \right \} =0 ,
\end{equation}

\noindent where the {\em inclinations} are defined as:

\begin{equation}
I_n^x \equiv{\mid v_n^x \mid \over \mid v_n^u \mid}, \quad I_n^s
\equiv{\mid v_n^s \mid \over \mid v_n^u \mid}, \label{inclin}
\end{equation}

\noindent
and

\begin{equation}
\label{stretch} {\mid v_{n+1}^u \mid \over \mid v_n^u \mid} >1.
\end{equation}

\noindent (\ref{straight}) implies that under iteration arbitrary
tangent vectors align with tangent vectors to the unstable
manifold and (\ref{stretch}) implies that these tangent vectors
also grow in length.

The proof essentially involves three steps. First we prove
(\ref{stretch}) for tangent vectors in $\Delta_n \cap W^s \left( M
\right)$. Next we extend this result to tangent vectors in
$\Delta_n$. Finally, we prove (\ref{stretch}).

\subsection{Inclinations for tangent vectors in the stable manifold}

We first prove (\ref{straight}) for $v_0$ in the tangent bundle of
$\Delta_n \cap W^s \left( M\right)$, denoted $T_{W^s \left(M
\right)} \Delta_n$. Let $p\in W^s (M)$, then by the invariance
property of the stable manifold (\ref{jacob1}) simplifies to:

\begin{equation}
Df_p = \left (
\begin{array}{lll}
A_s (x)+\partial_s r_s (p)& \partial_u r_s (p)& \partial_x r_s (p)+\partial_x A_s (x)\, s \\
0 & A_u (x ) +\partial_u r_u (p)& 0 \\
0 & \partial_u r_x (p)& \partial_x g (x)
\end{array}
\right ) . \label{jacob2}
\end{equation}

\noindent Acting on the tangent vector $(v_n^s, v_n^u, v_n^x)$, we
obtain the following relations

\begin{eqnarray}
v_{n+1}^s  & = & (A_s (x_n )+\partial_s r_s (p_n) ) v_n^s +\partial_u r_s (p_n )\, v_n^u + (\partial_x r_s (p_n)
+\partial_x A_s (x_n )\, s_n )v_n^x , \label{vsn} \\
v_{n+1}^u & = & (A_u (x_n )+\partial_u r_u (p_n ) ) v_n^u , \label{vun} \\
v_{n+1}^x & = & \partial_x r_u (p_n ) \,v_n^u + \partial_x g (x_n ) \, v_n^x .
\label{vxn}
\end{eqnarray}

\noindent
 Using these expressions, along with the estimates (\ref{k_est}), (\ref{lk1}), and (\ref{lk2}), we then obtain

\begin{eqnarray}
\mid v_{n+1}^s \mid  & \leq & (\lambda +k ) \mid v_n^s \mid +\parallel \partial_u
r_s (p_n) \parallel \, \mid v_n^u \mid + (k+ \parallel \partial_x A_s (x_n) \parallel \, \mid s_n\mid )
\mid v_n^x \mid , \label{vsn_est}\\
\mid v_{n+1}^u \mid & \geq & (\lambda^{-1} -k ) \mid v_n^u \mid , \label{vun_est}\\
\mid v_{n+1}^x \mid & \leq  & \parallel \partial_u r_x (p_n) \parallel \, \mid
v_n^u \mid + k \mid v_n^x \mid . \label{vxn_est}
\end{eqnarray}

\subsubsection{Inclination in the tangential direction}

Using (\ref{vxn_est}), (\ref{vun_est}) and (\ref{lk2}) gives:

\begin{equation}
{\mid v_{n+1}^x \mid \over \mid v_{n+1}^u \mid } \leq {k\over \lambda^{-1} -k}
{\mid v_n^x \mid \over \mid v_n^u \mid } +\parallel \partial_u r_x (p_n )\parallel .
\end{equation}

\noindent
 Using the estimate(\ref{secondorder}) on the second derivatives with the mean value
 inequality gives:

\begin{equation}
\label{lastx}
{\mid v_{n+1}^x \mid \over \mid v_{n+1}^u \mid } \leq {k\over
\lambda^{-1} -k} {\mid v_n^x \mid \over \mid v_n^u \mid } +C \mid s_n \mid .
\end{equation}

\noindent
 Let $p=(s,0,x)\in W^s (M)\cap U$ and $p_n =f^n (p)=(s_n
,0,x_n )$. By definition, we have $s_{n+1} =A_s (x_n ) s_n +r_s
(s_n ,0,x_n )$. Estimating this expression using assumption d) of
section \ref{normalform}, as well as the mean value inequality
with (\ref{k_est}) gives:

\begin{equation}
\mid s_{n+1} \mid \leq (\lambda +k) \mid s_n \mid ,
\end{equation}

\noindent from which it follows that:

\begin{equation}
\mid s_n \mid \leq (\lambda +k)^n \mid s\mid . \label{sit_est}
\end{equation}

\noindent Replacing $\mid s_n\mid$ by this expression in
(\ref{lastx}), we obtain

\begin{equation}
{\mid v_{n+1}^x \mid \over \mid v_{n+1}^u \mid } \leq {k\over \lambda^{-1} -k}
{\mid v_n^x \mid \over \mid v_n^u \mid } +C \mid s\mid (\lambda +k)^n .
\end{equation}

\noindent As a consequence, we have

\begin{equation}
{\mid v_n^x \mid \over \mid v_n^u \mid } \leq \left ( {k\over \lambda^{-1} -k}
\right ) ^n {\mid v_0^x \mid \over \mid v_0^u \mid } +C \mid s\mid
\sum_{i=0}^{n-1} \left ( {k\over \lambda^{-1} -k} \right ) ^i (\lambda +k)^{n-1
-i} .
\end{equation}

\noindent From (\ref{lk1}) and (\ref{lk2}) it follows that:

\begin{equation}
{k\over \lambda^{-1} -k} < k+\lambda ,
\end{equation}

\noindent Using this we obtain

\begin{equation}
{\mid v_n^x \mid \over \mid v_n^u \mid } \leq \left ( {k\over
\lambda^{-1} -k} \right ) ^n {\mid v_0^x \mid \over \mid v_0^u
\mid } +C \mid s\mid n (\lambda +k )^{n-1} .
\label{unstab_inc_est1}
\end{equation}

\noindent From (\ref{lk1})  we have $\lambda +k <1$, and
therefore:

\begin{equation}
\lim_{n\rightarrow +\infty} {\mid v_n^x \mid \over \mid v_n^u \mid } =0 .
\end{equation}

\subsubsection{Inclination in the stable direction}

Using (\ref{vsn_est}) and (\ref{vun_est}), we obtain

\begin{equation}
{\mid v_{n+1}^s \mid \over \mid v_{n+1}^u \mid} \leq \left ( {\lambda +k \over
\lambda^{-1}-k}\right ) {\mid v_n^s \mid \over \mid v_n^u \mid} +\parallel
\partial_u r_s (p_n )\parallel + (k+ \parallel \partial_x A_s (x_n )\parallel\, \mid s_n\mid ) {\mid v_n^x
\mid \over \mid v_n^u \mid} . \label{stab_inc_est1}
\end{equation}

\noindent Using the assumption (\ref{secondorder}), the mean value
inequality, and (\ref{sit_est}) gives:

\begin{equation}
\parallel \partial_u r_s (s_n ,0,x_n )\parallel \leq C\, \mid s_n \mid \leq
C\, \mid s \mid \, (\lambda +k)^{n} .
\end{equation}

\noindent Moreover, recall that from assumption iii) in the
statement of the $\lambda$-lemma we have:

\begin{equation}
\parallel \partial_x A_s (x_n )\parallel \leq D .
\end{equation}

\noindent Using these two estimates, (\ref{stab_inc_est1})
becomes:

\begin{equation}
{\mid v_{n+1}^s \mid \over \mid v_{n+1}^u \mid} \leq \left (
{\lambda +k \over \lambda^{-1}-k}\right ) {\mid v_n^s \mid \over
\mid v_n^u \mid} + C\mid s\mid\, (\lambda +k)^{n} + (k+D\mid
s\mid\, (\lambda +k)^{n} ) {\mid v_n^x \mid \over \mid v_n^u \mid}
. \label{stab_inc_est2}
\end{equation}

\noindent As a preliminary step to estimating
(\ref{stab_inc_est2}), we first estimate the third term on the
right-hand-side of (\ref{stab_inc_est2}) using
(\ref{unstab_inc_est1}):

\begin{equation}
(k+D\mid s\mid\, (\lambda +k)^{n} ) {\mid v_n^x \mid \over \mid
v_n^u \mid} \le (k+D\mid s\mid\, (\lambda +k)^{n} )\left(\left (
{k\over \lambda^{-1} -k} \right ) ^n {\mid v_0^x \mid \over \mid
v_0^u \mid } +C \mid s\mid n (\lambda +k )^{n-1}  \right)
\end{equation}

\noindent Now for $n$ sufficiently large we have:

\begin{equation}
k+D\mid s\mid\, (\lambda +k)^{n} <1,
\end{equation}

\noindent and by assumption (\ref{lk2}) we have $\lambda^{-1} -k
>1$, and therefore:

\begin{equation}
(k+D\mid s\mid\, (\lambda +k)^{n} ) {\mid v_n^x \mid \over \mid
v_n^u \mid} \le k^n {\mid v_0^x \mid \over \mid v_0^u \mid } +C
\mid s\mid n (\lambda +k )^{n-1}.
\end{equation}

\noindent Substituting this expression into (\ref{stab_inc_est2})
gives:

\begin{eqnarray}
{\mid v_{n+1}^s \mid \over \mid v_{n+1}^u \mid} & \leq & \left (
{\lambda +k \over \lambda^{-1}-k}\right ) {\mid v_n^s \mid \over
\mid v_n^u \mid} + C\mid s\mid\, (\lambda +k)^{n}  +C \mid s\mid n
(\lambda +k )^{n-1}+ k^n {\mid v_0^x \mid \over \mid v_0^u \mid },
\nonumber \\
 & \leq &  \left (
{\lambda +k \over \lambda^{-1}-k}\right ) {\mid v_n^s \mid \over
\mid v_n^u \mid} + \left(\lambda + k  \right)^{n-1} \left( C\mid
s\mid (1+ n) +  {\mid v_0^x \mid \over \mid v_0^u \mid } \right)
 \label{stab_inc_est3}
\end{eqnarray}

\noindent From this expression we obtain:

\begin{eqnarray}
{\mid v_{n+1}^s \mid \over \mid v_{n+1}^u \mid} & \leq & \left (
{\lambda +k \over \lambda^{-1}-k}\right )^n {\mid v_0^s \mid \over
\mid v_0^u \mid} + \sum_{i=0}^{n-1} \left ( {\lambda +k \over
\lambda^{-1}-k}\right )^i \left( \left(\lambda + k \right)^{n-i-2}
\left( C\mid s\mid (n-i) +  {\mid v_0^x \mid \over \mid v_0^u \mid
} \right)\right), \nonumber \\
 & \leq & \left (
{\lambda +k \over \lambda^{-1}-k}\right )^n {\mid v_0^s \mid \over
\mid v_0^u \mid} + \left(\lambda + k \right)^{n-2} n \left(C \mid
s \mid +  {\mid v_0^x \mid \over \mid v_0^u \mid }\right)
\end{eqnarray}

\noindent Consequently, as $n$ goes to infinity, we have

\begin{equation}
\lim_{n\rightarrow \infty} \di {\mid v_n^s \mid \over \mid v_n^u \mid } =0 .
\end{equation}

\subsection{Extending the estimates to tangent vectors not in $T_{W^s \left( M \right)} \Delta_n$}

We have shown that for $n$ sufficiently large, $U$ sufficiently
small,  for all $v\in T_{W^s \left( M \right)} \Delta_n$, $\mid
v\mid =1$,

\begin{equation}
I(v) \leq \epsilon .
\end{equation}

\noindent
 By continuity of the tangent plane, there exists
$\tilde{\Delta}_n \subset \Delta_n$ such that

\begin{equation}
\forall v\in T \tilde{\Delta}_n ,\ \  I(v)\leq 2\epsilon .
\end{equation}

\noindent For $n$ sufficiently large $\Delta_n$ is very close to
$W^u (M)$. We choose a neighbourhood $V_{\epsilon_s}$ of $W^u (M)$
of the form

\begin{equation}
V_{\epsilon_s} =\{ (s,u,x) \in B^s_{\epsilon_s} \times B^u_{\rho} \times M \} ,
\end{equation}

\noindent for $0<\epsilon_s <1$. Recall the estimates

\begin{equation}
\forall p\in  \subset U,\ \ \left .
\begin{array}{ll}
\parallel \partial_u r_x (p)\parallel \leq \delta ,& \ \parallel \partial_x r_x (p)\parallel \leq \delta,\\
\parallel \partial_u r_s (p)\parallel \leq \delta , & \ \parallel \partial_x r_s (p)\parallel \leq \delta
.
\end{array}
\right .
\end{equation}

\noindent Using  the mean value theorem and (\ref{secondorder}),
we choose $\epsilon_s$ such that:

\begin{equation}
\delta \geq C\epsilon_s .
\end{equation}

\noindent
 As a consequence, we can take

\begin{equation}
\delta =(C+1)\epsilon_s .
\end{equation}

\noindent Using (\ref{jacob1}), analogously to the estimates above
we obtain the following:

\begin{eqnarray}
\mid v_{n+1}^s \mid & \leq & (\lambda +k) \mid v_n^s\mid +\delta \mid v_n^u \mid +(C \epsilon_s +\delta ) \mid v_n^x \mid , \label{vsn_est2}\\
\mid v_{n+1}^u \mid & \geq & (\lambda^{-1} -k) \mid v_n^u \mid
-k\mid v_n^s \mid - (k+\rho C )\mid v_n^x \mid , \label{vun_est2}\\
\mid v_{n+1}^x \mid & \leq & k \mid v_n^s\mid +\delta \mid v_n^u
\mid +(\parallel \partial_x g (x_n )\parallel +\delta )\mid v_n^x \mid .
\label{vxn_est2}
\end{eqnarray}

\noindent We use these expressions to obtain the following
estimates of the inclinations:

\begin{eqnarray}
I_{n+1}^x & \leq & \di {1\over \lambda^{-1} -k} \left [ kI_n^s +\delta +(\parallel \partial_x g (x_n)\parallel
+\delta ) I_n^x \right ] \di \mu , \label{x_inc_est_4}\\
I_{n+1}^s & \leq & \di {1\over \lambda^{-1} -k} \left [ (\lambda
+k) I_n^s +\delta +(C\epsilon_s +\delta )I_n^x \right ] \di\mu ,
\label{s_inc_est_4}
\end{eqnarray}

\noindent where

\begin{equation}
\mu^{-1} = 1 - {k\over \lambda^{-1} -k} I_n^s -{k+C\rho \over
\lambda^{-1} -k} I_n^x .
\end{equation}

\noindent Since $\delta =(C+1)\epsilon_s$, and $\tilde{C}\geq \sup
\{
\parallel
\partial^2_{i,x} g (z)\parallel ,\ z\in U_{\rho} \}$, we obtain the
following estimates:

\begin{eqnarray}
(\parallel \partial_x g (x_n )\parallel +\delta ) I_n^x & \leq & \epsilon_s (\tilde{C}+C+1) I_n^x ,\\
\delta +(C\epsilon_s +\delta )I_n^x & \leq & \epsilon_s \left [
C+1 +I_n^x (2C+1) \right ].
\end{eqnarray}

\noindent We substitute these estimates into (\ref{x_inc_est_4})
and (\ref{s_inc_est_4}), and {\em assuming} that

\begin{equation}
I_n^x \leq \epsilon \quad \mbox{and} \quad I_n^s \leq \epsilon,
\label{inc_eps}
\end{equation}

\noindent we obtain:

\begin{eqnarray}
I^x_{n+1} & \le & \frac{\mu_*}{\lambda^{-1} -k} \left(k \epsilon +
(C+1) \epsilon_s + \epsilon \epsilon_s (M+C+1) \right),
\label{x_inc_est_5}
\\ I^s_{n+1} & \le & \frac{\mu_*}{\lambda^{-1} -k} \left(
(\lambda +k) \epsilon +  \epsilon_s  (C+1 +(2C+1) \epsilon
\right), \label{s_inc_est_5}
\end{eqnarray}

\noindent where $\mu_*^{-1} =1 - \di {2k+C\rho \over \lambda^{-1}
-k} \epsilon$.

\noindent Now if we choose $k$ small enough such that:

\begin{equation}
\frac{\lambda + k}{\lambda^{-1} -k} \mu_*<1,
\end{equation}

\noindent and $\epsilon_s$ satisfies:

\begin{equation}
\epsilon_s \le {\rm inf} \left\{
\begin{array}{l}
\epsilon \left(\frac{\lambda^{-1} -k}{\mu_*} \right)
\left(1-\frac{k \mu_*}{\lambda^{-1} -k} \right) \left(\frac{1}{C+1
+\epsilon(M+C+1)} \right), \\
\hfill \\
 \epsilon \left(1-\frac{\lambda +
k}{\lambda^{-1} -k}\mu_* \right) \frac{1}{C+1 +(2C+1) \epsilon}
\end{array}
\right.
\end{equation}

\noindent then

\begin{eqnarray}
I_{n+1}^x & \leq & \epsilon ,\\
I_{n+1}^s & \leq &  \epsilon .
\end{eqnarray}

\noindent Hence, we have shown that for $k$ and $\epsilon_s$
sufficiently small, $I_{n+1}^x \leq \epsilon$ and $I_{n+1}^s \leq
\epsilon$. Therefore the estimates
$I_n^s \leq \epsilon$ and $I_n^s \leq \epsilon$ are maintained under iteration.\\

As $\epsilon$ can be chosen as small as we want, the inclinations
$I_n^s$ and $I_n^x$ is as small as we want for $n$ sufficiently
large.

\subsection{Stretching along the unstable manifold}

We want to show that $f^n(\Delta ) \cap U$ is stretched in the
direction $W^u (M)$. In order to see this we compare the norm of a
tangent vector in $\tilde{\Delta}_n$ with its image under $Df$:

\begin{equation}
\sqrt{\frac{\mid v^s_{n+1} \mid^2 + \mid v^u_{n+1} \mid^2 + \mid
v^x_{n+1} \mid^2}{\mid v^s_{n} \mid^2 + \mid v^u_{n} \mid^2 + \mid
v^x_{n} \mid^2}} = \frac{\mid v^u_{n+1} \mid}{\mid v^u_{n} \mid}
\sqrt{ \frac{1+ \left(I^s_{n+1}\right)^2 +
\left(I^x_{n+1}\right)^2}{ 1+ \left(I^s_n \right)^2 + \left(I^x_n
\right)^2}}. \label{ratio}
\end{equation}

\noindent Using (\ref{vun_est2}) we obtain:

\begin{equation}
{\mid v_{n+1}^u \mid \over \mid v_n^u \mid} \geq \lambda^{-1} -k
-k\epsilon -(k+C\rho )\epsilon.
\end{equation}

\noindent Since $\epsilon$ can be chosen arbitrarily small we
have:

\begin{equation}
\lambda^{-1} -k -k\epsilon -(k+C\rho )\epsilon \geq \lambda^{-1} -2k >1 .
\end{equation}

\noindent Since the inclinations are arbitrarily small, it follows
that the norms of (nonzero) vectors in $\tilde{\Delta}_n$ are
growing by a ratio that approaches $\lambda^{-1} -2k >1$.
Therefore the diameter of $\tilde{\Delta}_n$ is increasing.
Putting this together with the fact that the tangent spaces have
uniformly small slope implies that there exists $\tilde{n}$ such
that for all $n \ge \tilde{n}$ $\tilde{\Delta}_n$ is $C^1$
$\epsilon$-close to $W^u(M)$.

This concludes the proof.

\end{document}